# A homogenized damping model for the propagation of elastic wave in a porous solid

Kangpei Meng, Q.M. Li*

Department of Mechanical, Aerospace and Civil Engineering, School of Engineering, the University of Manchester, Manchester, M13 9PL, UK

**Abstract:** This paper develops an averaging technique based on the combination of the eigenfunction expansion method and the collaboration method to investigate the multiple scattering effect of the SH wave propagation in a porous medium. The semi-analytical averaging technique is conducted using Monto Carlo method to understand the macroscopic dispersion and attenuation phenomena of the stress wave propagation in a porous solid caused by the multiple scattering effects. The averaging technique is verified by finite element analysis. Finally, a simple homogenized elastic model with damping is proposed to describe the macroscopic dispersion and attenuation effects of SH waves in porous media.

**Nomenclature**

| | |
|---|---|
| $a$ | The radius of the embedded cavities |
| $c$ | The wave speed in the matrix |
| $c_{\text{eff}}$ | The effective wave speed in the porous media |
| $k$ | The incident wavenumber |
| $k_{\text{eff}}$ | Effective wavenumber |
| $q$ | The number of the periodic mirror of cavities |
| $w^r$ | The displacement at a certain node |
| $\overline{w}^r$ | The sectional displacement averaged by displacements at nodes |
| $\overline{w}$ | The amplitude of sectional displacement |
| $\overline{w}^r_{\text{ran}}$ | The final averaged displacement at a certain section by the displacements at the sections in models with different distributions |
| $\overline{w}_{\text{ran}}$ | The amplitude of final averaged displacement at a certain section |
| $w_{\text{in}}$ | The amplitude of the incident wave |
| $E$ | The Young's modulus of the matrix |
| $E_{\text{eff}}$ | The effective Young's modulus of the porous media |
| $F$ | The load applied on nodes on the boundary between the established finite element model and the infinite elements |
| $H$ | The height of the selected representative segments |
| $H_m^{(1)}$ | The $m$th order of Hankel function of the kind |
| $J_m()$ | The $m$th order of Bessel function |
| $L$ | The number of random distributions |
| $M$ | The maximum order of Hankel function in calculation |
| $N$ | The total number of the embedded cavities |
| $Q$ | The repeating number of the representative element |
| $R_j$ | The distance of the selected point in the polar coordination system centred at $j$th cavity |
| $T$ | The length of the selected representative segments |
| $\theta_j$ | The angle of the selected point in the polar coordination system centred at $j$th cavity |
| $\lambda$ | The incident wavelength |
| $\lambda_{\text{eff}}$ | Effective wavelength in the porous media |
| $\mu$ | The shear modulus of the matrix |
| $\mu_{eff}$ | The effective shear modulus of the matrix |
| $\rho$ | The density of the matrix |
| $\rho_{\text{eff}}$ | Effective wavelength of the porous media |
| $\omega$ | The frequency of incident wave |

# 1. Introduction

Heterogeneous porous materials have attracted great attentions due to their lightweight and enhanced mechanical properties (e.g. high specific stiffness, strength and toughness) and various preferred functions (e.g. thermal and acoustic insulation). Porous materials are distinguished from condensed materials mainly due to the existence of mesoscale pores. In order to understand their mechanical behaviours, it is necessary to understand the effects of their mesoscale structures on the macroscopic mechanical responses, which is especially important for dynamic loading problems because stress wave propagation may interact with the mesoscale features in a porous material.

Stress and strain are valid physical quantities in all meso-scale phases of a heterogeneous medium. When a stress wave with high-frequency components propagates in a heterogeneous medium, multiple reflections of the stress wave happen on the boundaries of the inclusions (or pores for a porous material) in the heterogeneous material, leading to macroscopic dispersion and attenuation (Foldy, 1945; Lax, 1951, 1952). The dispersion effect denotes the change of wave speed in a medium with the frequency of the incident wave while the attenuation effect relates to the gradual decrease of wave amplitude during wave propagation. In general, the inclusions in a heterogeneous material have different sizes and irregular shapes. When the shape of the inclusion is cylindrical, spherical or ellipsoidal, multiple scattering method (MSM), introduced initially by Foldy (1945), can be approximately used to consider multiple wave reflections by the boundary of inclusions and evaluate the elastic wave propagation in a heterogeneous medium with discrete inclusions.

MSM (Foldy, 1945) was extended to more general applications, i.e. anisotropic scattering, inelastic scattering, and scattering of quantized waves including photons etc. and combined with 'quasi-crystalline approximation' by Lax (Lax, 1951, 1952) (thus known as Foldy-Lax method). Twersky (Twersky, 1952; Twersky, 1962a, b), Waterman and Truell (Waterman and Truell, 1961), Fikioris and Waterman (Fikioris and Waterman, 1964) and Varadan et al. (Varadan et al., 1985a; Varadan et al., 1985b; Varadan et al., 1978) made further contributions to the development of T-matrix method to obtain averaged wave propagation properties. Apart from T-matrix method, Bose and Mal (Bose and Mal, 1973, 1974) applied eigenfunction expansion method, i.e. using the Hankel functions of the first kind to express the wave scattering from the cylindrical scatters and transform all calculations into a single coordination system using Graf's addition theorem (Abramowitz and Stegun, 1964) to calculate the wave field in a 2D model with cylindrical inclusions. They then combined 'quasi-crystalline approximation' and their own 'pair correlation' (i.e. the probability densities decline in an exponential way) to obtain effective field quantities in a model with randomly distributed scatters. When scatters are distributed randomly in the matrix, for example dusts in the atmosphere, most researches calculate an average wave property based on the assumptions of 'quasi-crystalline approximation' and 'pair correlation'. However, the results from (Kim, 2010) showed that the choice of 'pair correlation' may impose considerable effect on the effective wavenumber during homogenisation, especially when the incident wave frequency is relatively high and the model has dense scatters.

Both T-matrix method and eigenfunction expansion method can give results with enough accuracy when the positions of embedded scatters are identified. But the computation becomes cumbersome for high incident wave frequency and dense scatter distribution. Biwa (Biwa et al., 2007; Biwa et al., 2004) and Sumiya et al (Sumiya et al., 2013) presented a computational procedure based on a collaboration method to determine the coefficients in T-matrix and gave a semi-analytical result of SH wave field in matrix with long cylindrical fibre, in which the whole-scale randomly-distributed scatters were replaced by periodic arrangements in the direction perpendicular to the wave propagation direction. The wave field of the whole model can be treated as periodic if the number of periodic arrangements is large enough, which significantly reduce the computational load. The combination of partial periodicity and Monte-Carlo study facilitates the calculation of the averaged wave propagation properties of the model with random scatters.

MSM has been shown as a useful tool to study the macroscopic wave dispersion and attenuation phenomena originated from the mesoscale heterogeneity of the material. MSM's applications have been extended to the wave propagations in heterogeneous media with considering more complex inclusion-matrix systems, e.g. the imperfect bonding (Zhang et al., 2018), coated fibre (Biwa and Sumiya, 2015) and the transition zone (Ghanei Mohammadi and Hosseini Tehrani, 2019) in a composite. It is desirable that the effects of the inclusions in a heterogeneous material could be represented by the effective parameters of its corresponding homogenized counterpart material (e.g. using effective modulus, Poisson's ratio, density etc.) and other pseudo effective parameters (e.g. damping). Therefore, the macroscopic stress wave propagation could be described by a homogenized constitutive model without the need of mesoscale modelling of the heterogeneous material, which can significantly reduce the computational time. However, it is still a challenge to homogenize material responses based on mesoscale modelling (e.g. FE) due to its high demand of computational resources. MSM offers an approximate, but more efficient, approach to homogenize heterogonous materials.

This paper employs MSM, i.e. combining the eigenfunction expansion method adopted by Bose (Bose and Mal, 1973) and the collaboration-method-based procedure proposed by Biwa (Biwa et al., 2004), to investigate the SH wave propagation in a heterogeneous medium. Particularly, porous medium is considered for computational convenience without losing the generality of the methodology. The semi-analytical MSM method is verified using finite element model (FEM) Abaqus and an averaging process based on a Monto Carlo study. Finally, a homogenized elastic model with structural damping is proposed to describe the macroscopic dispersion and attenuation of SH waves in porous materials.

2. **Methodology**

2.1 Multiple scattering problem

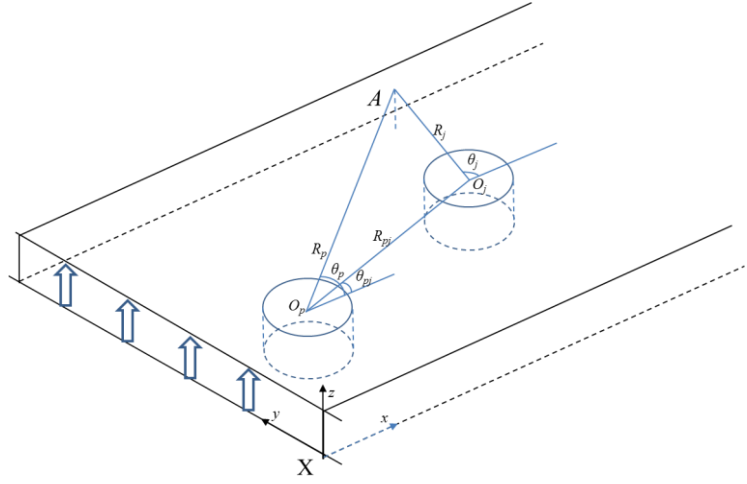

Fig. 1. Multiple scattering coordinate system.

Consider the case that there are more than one cylindrical cavities embedded in an infinite matrix material in coordinate system *X*, as shown in Fig. 1, where a time-harmonic SH plane wave propagates in *x* direction. The total displacement at a specific point in the material domain can be obtained as the superposition of the incident and scattered displacement waves, at the same point, from all cavity surfaces.

Suppose the shear modulus and the density of the matrix are respectively $\mu$ and $\rho$, and all cavities have the same radius *a*. Taking a point A outside the cavities, the total displacement after suppressing the time factor $e^{-i\omega t}$ can be expressed as (Bose and Mal, 1973)

$$w = e^{ikx} + \sum_{j=1}^{N} \sum_{m=-\infty}^{\infty} A_{jm} H_m^{(1)}(kR_j) e^{im\theta_j} \qquad (1)$$

where *k* is the incident shear wavenumber in the matrix; *x* is the position coordinate of point A in the global coordination system *X* in Fig. 1. The first part on the right side of Eq.(1) is the incident wave with unit amplitude, and the second part represents the scattered waves from the cavities. *N* is the total number of cavities. $A_{jm}$ is a series of constants to be determined. $H_m^{(1)}()$ is the Hankel function of the first kind in *m*th order. $R_j$ is the distance between the centre of *j*th cavity and point A. $\theta_j$ is the angle between A in *j*th cavity-centred coordinate system and the wave propagation direction. ($R_j$, $\theta_j$) denotes the position of point A in *j*th cavity-centred coordinate system.

For the convenience of further calculation, Graf's addition theorem (Abramowitz and Stegun, 1964) is employed to transform all calculations for point A from the *p*th cavity-centred coordinate system to the *j*th cavity-centred coordinate system, i.e.

$$H_n^{(1)}(kR_p) e^{in\theta_p} = \sum_{m=-\infty}^{\infty} (-1)^m H_{m+n}^{(1)}(kR_{jp}) J_m(kR_j) e^{i(m+n)\theta_{jp}} e^{-im\theta_j} \qquad (2)$$

where *n* is an arbitrary integer number; ($R_p, \theta_p$) is the coordinate of point A in *p*th cavity-centred coordinate system; ($R_{jp}, \theta_{jp}$) is the coordinate of *p*th cavity in the *j*th cavity-centred coordinate system. $J_m()$ is the Bessel function in *m*th order. Thus, all the coefficients of the Hankel function at *p*th cavity in Eq.(1) can be calculated to get the displacement of a point in *j*th cavity-centred coordinate system through the transformation from *p*th cavity-centred coordinate system to *j*th cavity-centred coordinate system.

The study of wave scattering in porous media can be treated as a special case of the model with fibre inclusions. The stress fields can be obtained from geometrical equations and Hooke's law, which have been given in (Bose and Mal, 1973) and (Mow and Pao, 1971), and therefore, will not be presented here. The zero shear stress condition on the surface of *j*th cavity gives

$$0 = \mu \frac{\partial w}{\partial R_j}|_{R_j=a} = \left\{ \mu \cdot A_{jn} \frac{\partial}{\partial R_j} H_n(kR_j) + \mu \cdot \frac{\partial}{\partial R_j} J_n(kR_j)[i^n e^{ikR_j \cos(\theta_j)} \right. \\ \left. + \sum_{p=1,p\neq j}^{N} \sum_{m=-\infty}^{\infty} (-1)^m H_{m+n}^{(1)}(kR_{jp}) e^{i(m+n)\theta_{jp}} e^{-im\theta_j}] \right\}|_{R_j=a} \tag{3}$$

Combining Eqs. (2) and (3), $A_{jn}$ can be calculated, i.e.

$$A_{jn} = iB_n C_{jn} \tag{4}$$

where

$$iB_n = -\frac{\frac{\partial}{\partial a} J_n(ka)}{\frac{\partial}{\partial a} H_n(ka)} \tag{5}$$

and

$$C_{jn} = i^n e^{ikr_j \cos\theta_j} + i \sum_{p=1,p\neq j}^{N} \sum_{m=-\infty}^{\infty} B_m C_{pm} H_{m-n}(kR_{pj}) e^{i(m-n)\theta_{pj}} \tag{6}$$

Eq. (6) can be transformed into a set of linear equations after selecting approximate truncation, which can then be solved. The total displacement at A can be calculated by combining Eq. (1) and Eqs.(4-6).

2.2 Homogenization of the dynamic behaviour of a porous medium with random cavities

Wave attenuation has been observed when a time-harmonic wave propagates in heterogeneous materials (Aggelis and Philippidis, 2004; Biwa et al., 2004; Bose and Mal, 1973; Huang and Rokhlin, 1995; Philippidis and Aggelis, 2003; Tauchert and Guzelsu, 1972). To quantify the wave attenuation in a porous medium, an elastic matrix model embedded with randomly-distributed cylindrical cavities is established. It is extremely time-consuming to use MSM to calculate the wave propagation in a sufficiently large-scale material domain with randomly-distributed cavities, which will be considered in this study and is schematically shown in Fig. 2 (a). Thus, a quasi-periodic description, similar to the method used in (Biwa et al., 2004), is adopted in the present study to reduce the computational demands. A rectangular segment with length 40 mm and width 20 mm is taken as a representative segment in the in-plane coordinate system *X* (*x*, *y*), as shown in Fig. 2 (b). Set the wave propagation direction as positive *x* direction and establish a coordinate system with the origin at the middle point of the left boundary of the representative segment. 50 cavities are distributed randomly in the representative segment matrix. The rectangular representative segment is repeated 300 times in both positive and negative *y* directions (i.e. there are totally 601 repeated representative segments). The radius of the cavities keeps constant of 0.6 mm. The shear modulus and density of the matrix are respectively set as 26.92 GPa and $2.7\times10^3$ kg/m$^3$.

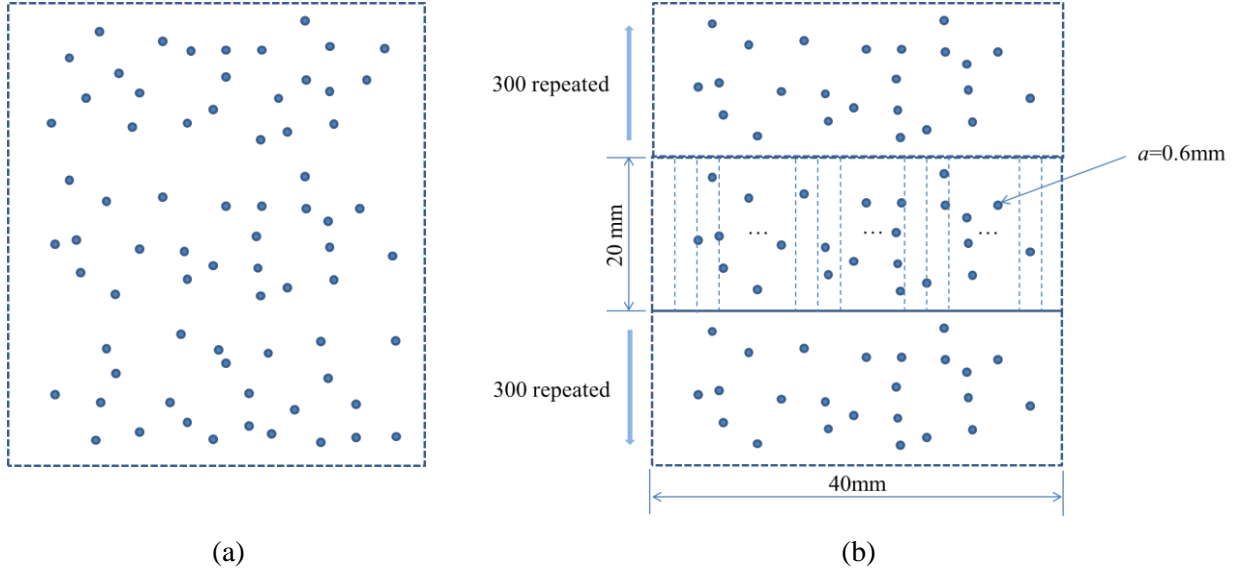

(a)              (b)

Fig. 2. A schematic drawing for multiple scattering model (a) random distribution; (b) periodic random distribution.

Then, the coefficients $C_{jn}$ in Eq. (6) for the model in Fig. 2 (b) can be obtained, i.e.

$$C_{jn} = i^n e^{ikr_j \cos\theta_j} + i \sum_{p=1,p\neq j}^{N} \sum_{q=-Q}^{Q} \sum_{m=-M}^{M} B_m C_{pm} H_{m-n}(kR_{pqj}) e^{i(m-n)\theta_{pqj}}$$
$$+ i \sum_{q=-Q,q\neq 0}^{Q} \sum_{m=-M}^{M} B_m C_{jm} H_{m-n}(kR_{jqj}) e^{i(m-n)\theta_{jqj}} \tag{7}$$

where

$$R_{pqj} = \sqrt{(x_p - x_j)^2 + (y_p - y_j + q \times H)^2} \tag{8}$$

$$\theta_{pqj} = \arccos \frac{x_p - x_j}{\sqrt{(x_p - x_j)^2 + (y_p - y_j + q \times H)^2}} \tag{9}$$

$$R_{jqj} = |q \times H| \tag{10}$$

and

$$\theta_{jqj} = \arccos \frac{x_j - x_j}{\sqrt{(x_j - x_j)^2 + (y_j - y_j + q \times H)^2}} = \frac{\pi}{2} \tag{11}$$

in which, $R_{pqj}$ is the distance between the *p*th cavity and the *q*th periodic mirror of *j*th cavity; $\theta_{pqj}$ is the angle between the *p*th cavity and the *q*th periodic mirror of *j*th cavity; $R_{jqj}$ is the distance between the *j*th cavity and the *q*th periodic mirror of *j*th cavity; $\theta_{jqj}$ is the angle between the *j*th cavity and the $q$th periodic mirror of $j$th cavity; $H$ and $T$ are the width and length of the rectangular representative segment, which are taken as 20mm and 40mm, respectively, in this study; $Q$ is the repeating number of the rectangular representative segment, which is taken as 300; $M$ is the maximum order of Hankel functions in Eq. (7), which is taken as 10 based on a convergence study; $N$ is the total number of cavities in the representative segment, which is taken as 50. Eq. (7) can be transformed into a set of linear equations after truncation and then be solved.

Combining Eqs. (1), (4), (5) and (7), the actual dynamic response of every point in the material domain could be calculated. The calculation result at every point is a complex number in the form of $w^r + w^i i$, from which both the displacement and phase angle at a given point can be determined. The amplitude of the displacement at that point is $\bar{w} = \sqrt{(w^r)^2 + (w^i)^2}$.

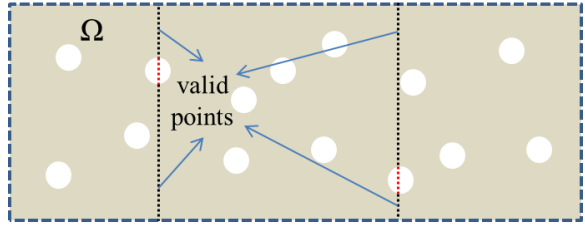

Fig. 3. A schematic drawing for selected points to average the displacement.

Take set $\Omega$ as the domain space outside the cavities in the central representative segment shown in Fig. 3. The $z$-directional displacements at 400 equally-spaced $x$ coordinate sections and 100 equally-interval points in $y$ direction at the same $x$ coordinate section in the $\Omega$-domain of the central representative segment (only the points outside the cavities are used to obtain the displacement) are calculated to study the wave dispersion and attenuation, i.e.

$$w^r(x_g, y_l) = \begin{cases} w^r(x_g, y_l), & (x_g, y_l) \in \Omega \\ 0, & (x_g, y_l) \notin \Omega \end{cases} \quad (12a)$$

$$w^i(x_g, y_l) = \begin{cases} w^i(x_g, y_l), & (x_g, y_l) \in \Omega \\ 0, & (x_g, y_l) \notin \Omega \end{cases} \quad (12b)$$

where $x_g = g \times \Delta x, y_l = l \times \Delta y$, $\Delta x = \frac{40\ mm}{400} = 0.1\ mm, \Delta y = \frac{20\ mm}{100} = 0.2\ mm$; $g$ ranges from 1 to 400 and $l$ from 1 to 100. The z-directional homogenized displacement at a given $x$ coordinate section (i.e. the sectional displacement components) is then calculated based on the calculation results of 100 equally-interval points in $y$ direction in the central representative segment at each of 400 $x$-coordinates, i.e.

$$\bar{w}^r(x_g) = \frac{1}{NUM(x_g)} \sum_{l=1}^{NUM} w^r(x_g, y_l) \quad (13a)$$

$$\bar{w}^i(x_g) = \frac{1}{NUM(x_g)} \sum_{l=1}^{NUM} w^i(x_g, y_l) \quad (13b)$$

where $NUM(x_g)$=100 is the number of the $y_l$ points when $x = x_g$. The variation of either real or imaginary component, i.e. $\bar{w}^r(x_g)$ or $\bar{w}^i(x_g)$, with all $x_g$ coordinates when g varies from 1 to 400 can be used to determine the effective wavelength of the sectional displacement wave. Therefore, the dependence of the phase velocity of the sectional displacement wave on the incident wavelength (or frequency) can be evaluated to study the wave dispersion due to the existence of cavities. Meanwhile, the amplitude of the sectional displacement at $x = x_g$ can be obtained from

$$\bar{w}(x_g) = \sqrt{[\bar{w}^r(x_g)]^2 + [\bar{w}^i(x_g)]^2} \tag{14}$$

whose variation with $x_g$ (g=1, 2, … …, 400) can be used to study the attenuation of the sectional displacement wave.

In order to have a better representation of the random distribution of cavities, the distributions of the cavities are generated *L* times randomly by the software MATLAB using random sequential adsorption (RSA) method (Feder, 1980). The calculations in Eq. (13a, b) are repeated *L* times (normally, repeating number of 20 is sufficient, i.e. *L*≤20) and averaged again to obtain final averaged displacement at $x = x_g$ and its amplitude in such periodic random distribution model, i.e.

$$\bar{w}^r_{\text{ran}}(x_g) = \frac{1}{L}[\bar{w}^r(x_g)(case1) + \bar{w}^r(x_g)(case2) + \cdots + \bar{w}^r(x_g)(caseL)] \tag{15a}$$

$$\bar{w}^i_{\text{ran}}(x_g) = \frac{1}{L}[\bar{w}^i(x_g)(case1) + \bar{w}^i(x_g)(case2) + \cdots + \bar{w}^i(x_g)(caseL)] \tag{15b}$$

$$\bar{w}_{\text{ran}}(x_g) = \sqrt{[\bar{w}^r_{\text{ran}}(x_g)]^2 + [\bar{w}^i_{\text{ran}}(x_g)]^2} \tag{16}$$

Similarly, the real (or imaginary) component and the amplitude of the final averaged displacement can be used to study the wave dispersion and attenuation, respectively, caused by the existence of cavities. It is well-known that wave dispersion and attenuation can be significantly influenced by the wavenumber of the incident wave. In this study, the wavenumber of the harmonic incident wave *k* will vary from 400/m to 2000/m, which corresponds to the dimensionless incident wavenumber *ka* from 0.24 to 1.20 and dimensionless incident wavelength, namely the ratio of the incident wavelength to the cavity diameter ($\lambda/(2a)$), from 13.10 to 2.61.

### 3. Macroscopic observations

3.1 Dispersion effect

The dispersion effect can be evaluated by comparing the incident wavenumber with the effective wavenumber in a porous medium where the effective wavenumber is calculated from its relationship with effective wavelength. The determination of the effective wavelengths will be discussed later in this sub-section. The (real component) curves of the dimensionless sectional displacement along *x* coordinate are calculated at 400 *x*-coordinates using Eq. (13a). Fig. 4 shows the dimensionless sectional displacement curves for 5 random cavity distributions generated by RSA method (i.e. *L*=5) when the dimensionless incident wavenumber *ka* increases from 0.24 to 1.2, where $\bar{w}^r$ is the real component of the sectional displacement and $w_{\text{in}}$ is the amplitude of the incident wave displacement, which is selected as unity in this study.

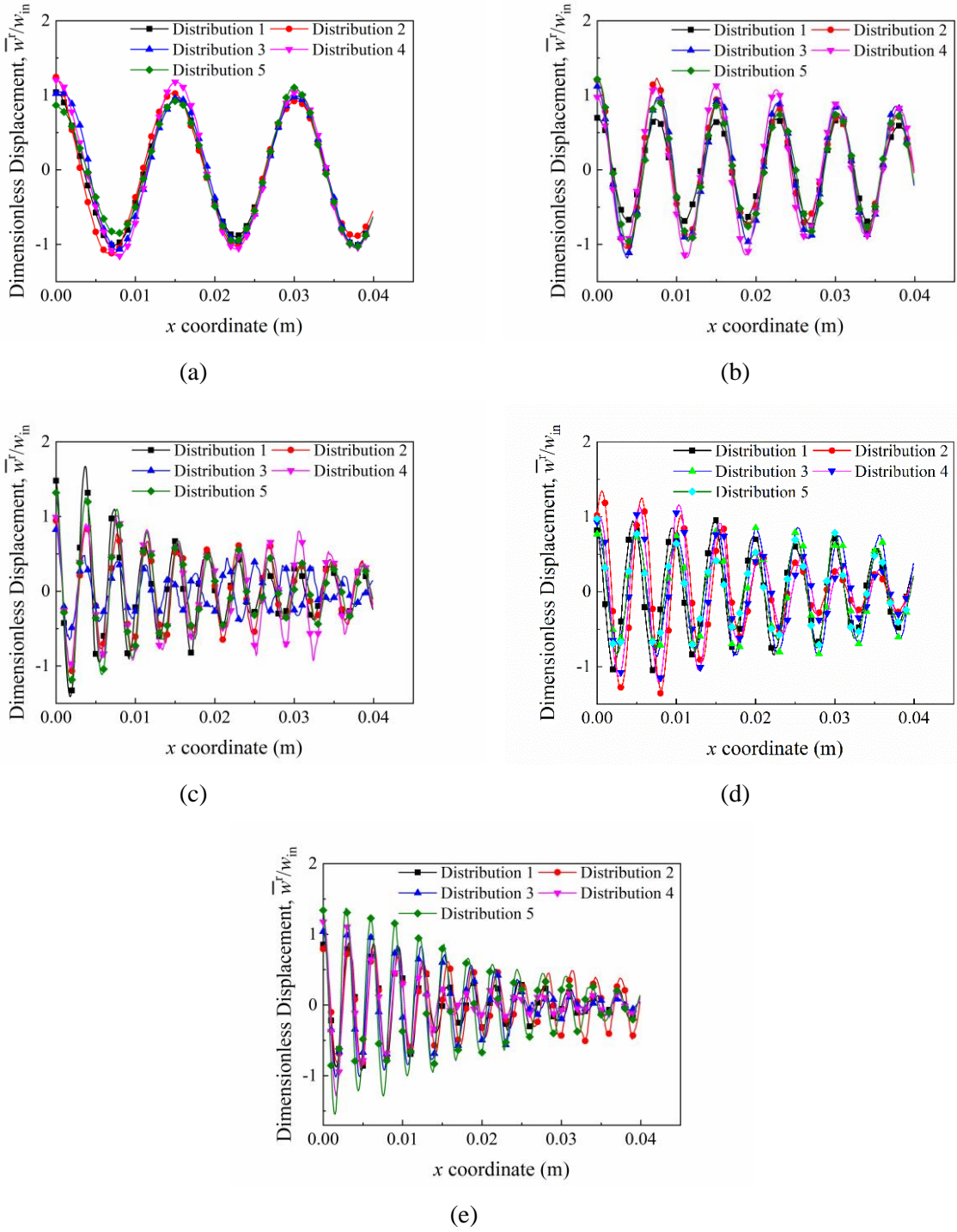

Fig. 4. The dimensionless sectional displacement (real component $\bar{w}^{\text{r}}/w_{\text{in}}$) along $x$ coordinate for each of the 5 random cavity distributions ($L$=5) and different dimensionless incident wavenumbers of (a) $ka$=0.24; (b) $ka$ =0.48; (c) $ka$ =0.72; (d) $ka$ =0.96; (e) $ka$ =1.20.

As shown in Fig. 4, when the dimensionless incident wavenumber $ka$=0.24, the sectional displacement curve along $x$ coordinate is insensitive to the cavity distributions. However, with the increase of the incident wavenumber, the ratios of the maximum of the first peak to the minimum of the first peak of the five curves are 1.44, 1.74, 1.83, 2.16 and 1.90 for *ka*=0.24, *ka*=0.48, *ka*=0.72, *ka*=0.96 and *ka*=1.20, respectively, which means that the sectional displacement curve along *x* coordinate generally becomes less dependent of the cavity distributions with the reduction of dimensionless incident wavenumber. When the dimensionless incident wavenumber *ka* approaches to zero, the ratio of the maximum of the first peak to the minimum of the first peak approaches to unity.

Fig. 5 shows further averaged curves among those of *L* random cavity distributions for *L*=5 and *L*=20, according to Eq. (15a).

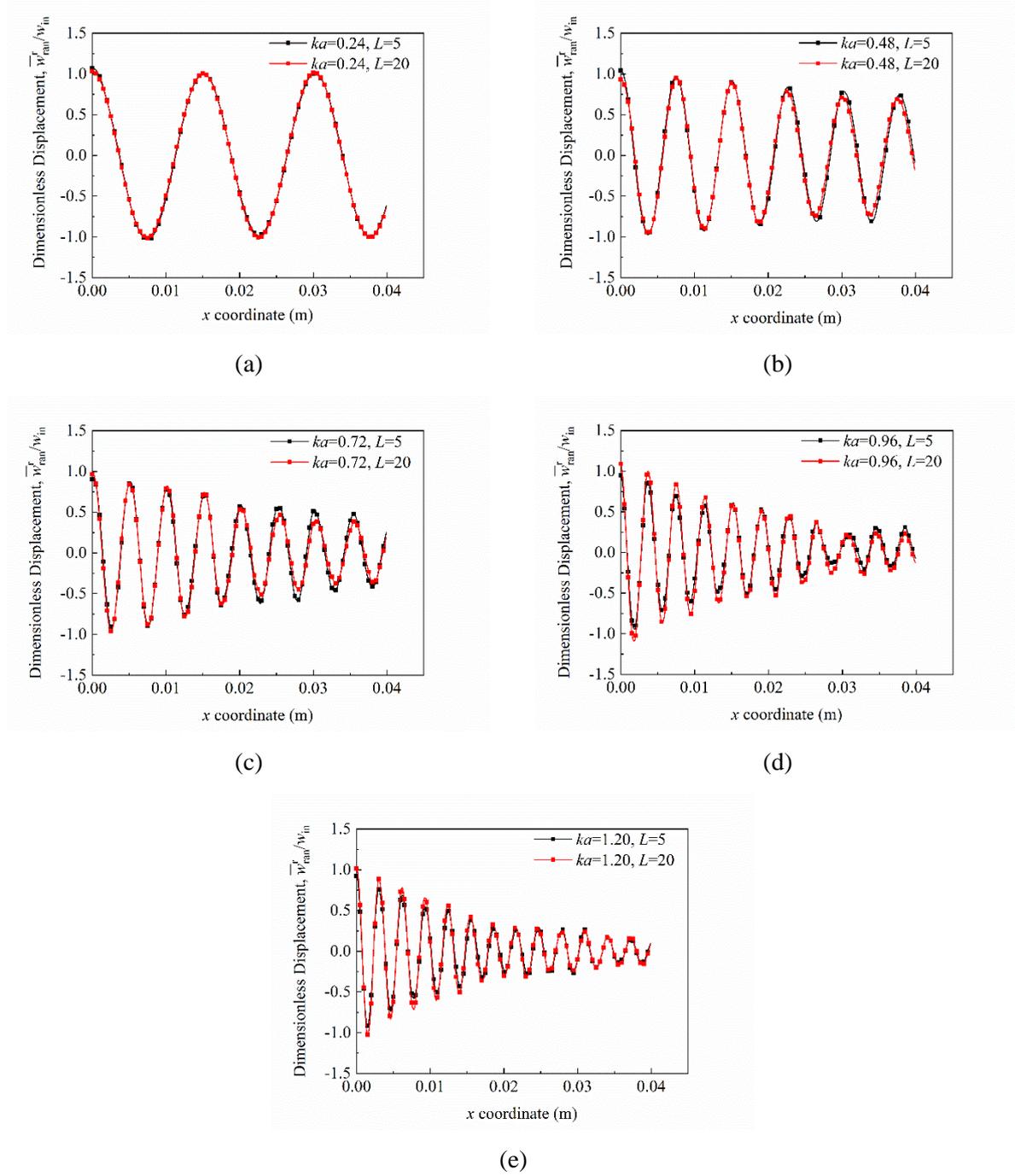

Fig. 5. The dimensionless final averaged displacement (real component $\bar{w}^r/w_{\text{in}}$) along *x* coordinates for *L*=5 and *L*=20, and different incident wavenumbers, (*a*) *ka*=0.24; (b) *ka*=0.48; (c) *ka*=0.72; (d) *ka*=0.96; (e) *ka*=1.20.

As shown in Fig. 5, when the dimensionless incident wavenumber $ka$=0.24, the dimensionless final averaged displacement wave has approximately a harmonic distribution along $x$ coordinate, which is independent of $L$ (i.e. they are nearly the same for $L$=5 and $L$=20). However, their differences gradually increase with the increase of $ka$, indicating the need of more random cavity distributions (i.e. larger L number) for higher values of $ka$ (i.e. relatively larger cavity size). On the other hand, with the increase of the dimensionless incident wavenumber toward $ka$=1.20, the averaged wavelength gradually decreases and the displacement amplitude gradually shows a more significant attenuation trend during propagation, which will be discussed further in Section 3.2.

The effective wavelength can be directly measured as the distance between the two nearby peaks on the dimensionless final averaged displacement curve against $x$ coordinate for different incident wavenumbers. A typical example for $ka$=1.20 is shown in Fig. 6.

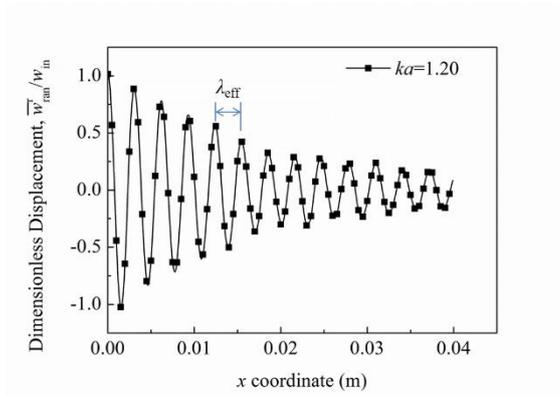

Fig. 6. The measurement of effetive wavelength for $ka$=1.20.

The effective wavenumber ($k_{\text{eff}}$) is related to the effective wavelength ($\lambda_{\text{eff}}$) by

$$k_{\text{eff}} = \frac{2\pi}{\lambda_{\text{eff}}} \tag{17}$$

The calculation results of effective wavelengths, effective wavenumbers and dimensionless effective wavenumbers for different incident wavenumbers are compared in Table 1.

Table 1. The comparisons between effective wavenumbers and incident wavenumbers.

| Dimensionless incident wavenumber $ka$ | Incident wavenumber $k$ | Effective wavelength $\lambda_{\text{eff}}$ | Effective wavenumber $k_{\text{eff}}$ | Dimensionless effective wavenumber $k_{\text{eff}}/k$ | Error to average $k_{\text{eff}}/k$ (1.02) |
|---|---|---|---|---|---|
| 0.24 | 400 /m | 0.0151 m | 416.1 /m | 1.04 | 2.0% |
| 0.48 | 800 /m | 0.0076 m | 826.7 /m | 1.03 | 1.0% |
| 0.72 | 1200 /m | 0.0051 m | 1231.0 /m | 1.03 | 1.0% |
| 0.96 | 1600 /m | 0.0039 m | 1611.1 /m | 1.01 | 1.0% |
| 1.20 | 2000 /m | 0.0031 m | 2026.8 /m | 1.01 | 1.0% |

As shown in Table 1, the values of the dimensionless effective wavenumber $k_{\text{eff}}/k$ are close to unity and the maximum error is 4% when $k_{\text{eff}}/k$ ranges from 0.24 to 1.20. The average value of $k_{\text{eff}}/k$ is 1.02 and the maximum error of $k_{\text{eff}}/k$ to the average value is less than 2%.

For an incident wave with circular frequency $\omega$, we have

$$c = \lambda\omega = \frac{\omega}{k} \qquad (18)$$

where $c$ is the phase speed and $\lambda$ is the wavelength. When the incident wave propagates in a porous medium, Eq.(18) is applicable to its phase speed, wavelength and wavenumber represented by their respective effective counterparts, i.e. $c_{\text{eff}}$, $\lambda_{\text{eff}}$ and $k_{\text{eff}}$, which can be calculated using MSM or meso-scale FEM.

According to Eq.(18) and Table 1, the effective phase speed in the studied porous medium is almost independent of the incident wavelength (or frequency), and therefore, there is negligible dispersion effect caused by the cavities in the porous medium, when the dimensionless incident wavenumber is less than 1.20 in the present case studied.

3.2 Attenuation effect

The attenuation effect is evaluated by fitting the declining amplitude of the final averaged displacement curve along $x$ coordinate.

Again, in line with Fig. 4 in Section 3.1, the dimensionless amplitude of the sectional displacement curves ($\bar{w}/w_{\text{in}}$) calculated from Eq. (14) at 400 $x$-coordinates (note: the symbols are shown with the interval of 10 $x$-coordinate) are shown in Fig. 7 for each of the 5 random cavity distributions ($L$=5) and different dimensionless incident wavenumbers $ka$ from 0.24 to 1.20.

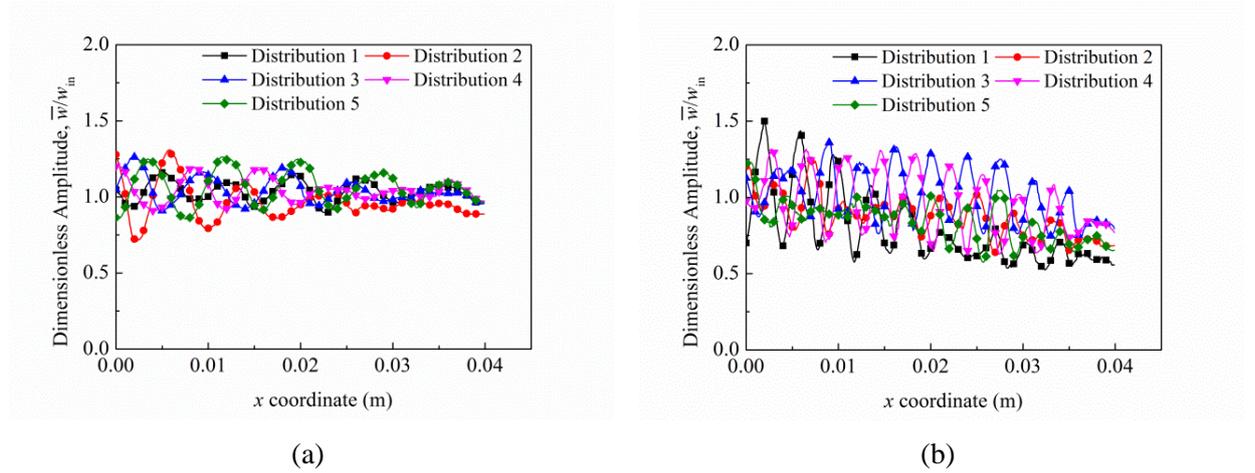

(a)                                    (b)

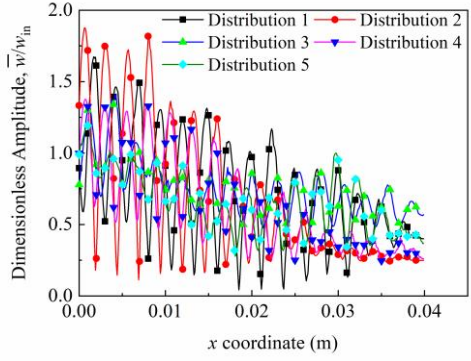
(c)

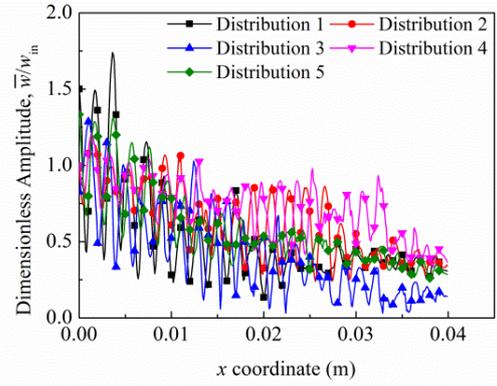
(d)

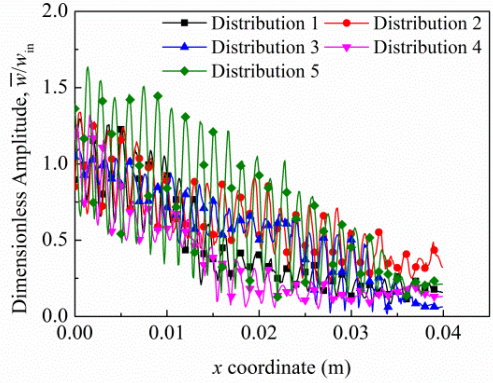
(e)

Fig. 7. The dimensionless amplitude of the sectional displacement ($\bar{w}/w_{\text{in}}$) along $x$ coordinates for each of the 5 random cavity distributions ($L$=5) and different dimensionless incident wavenumbers, (*a*) $ka$=0.24; (b) $ka$=0.48; (c) $ka$=0.72; (d) $ka$=0.96; (e) $ka$=1.20.

Similar to the results obtained in Section 3.1, a good repeatability in the dimensionless amplitude of the sectional displacement curves along $x$ coordinate is generally kept when the dimensionless incident wavenumber $ka$=0.24 and the difference among five curves generally becomes more significant with the increase of the dimensionless incident wavenumber. The ratios of the maximum value of the first peak to the minimum value of the first peak of five curves are 1.12, 1.25, 1.34, 1.24 and 1.49 at $ka$=0.24, $ka$=0.48, $ka$=0.72, $ka$=0.96 and $ka$=1.20, respectively, which shows that the difference among five curves with different cavity distributions generally increases with the incident wavenumber.

The dependence of the dimensionless amplitude of the final averaged displacement calculated by Eq.(16) on the number of the random cavity distributions ($L$) is shown in Fig. 8, where $L$ was taken as 5, 10, 15 and 20.

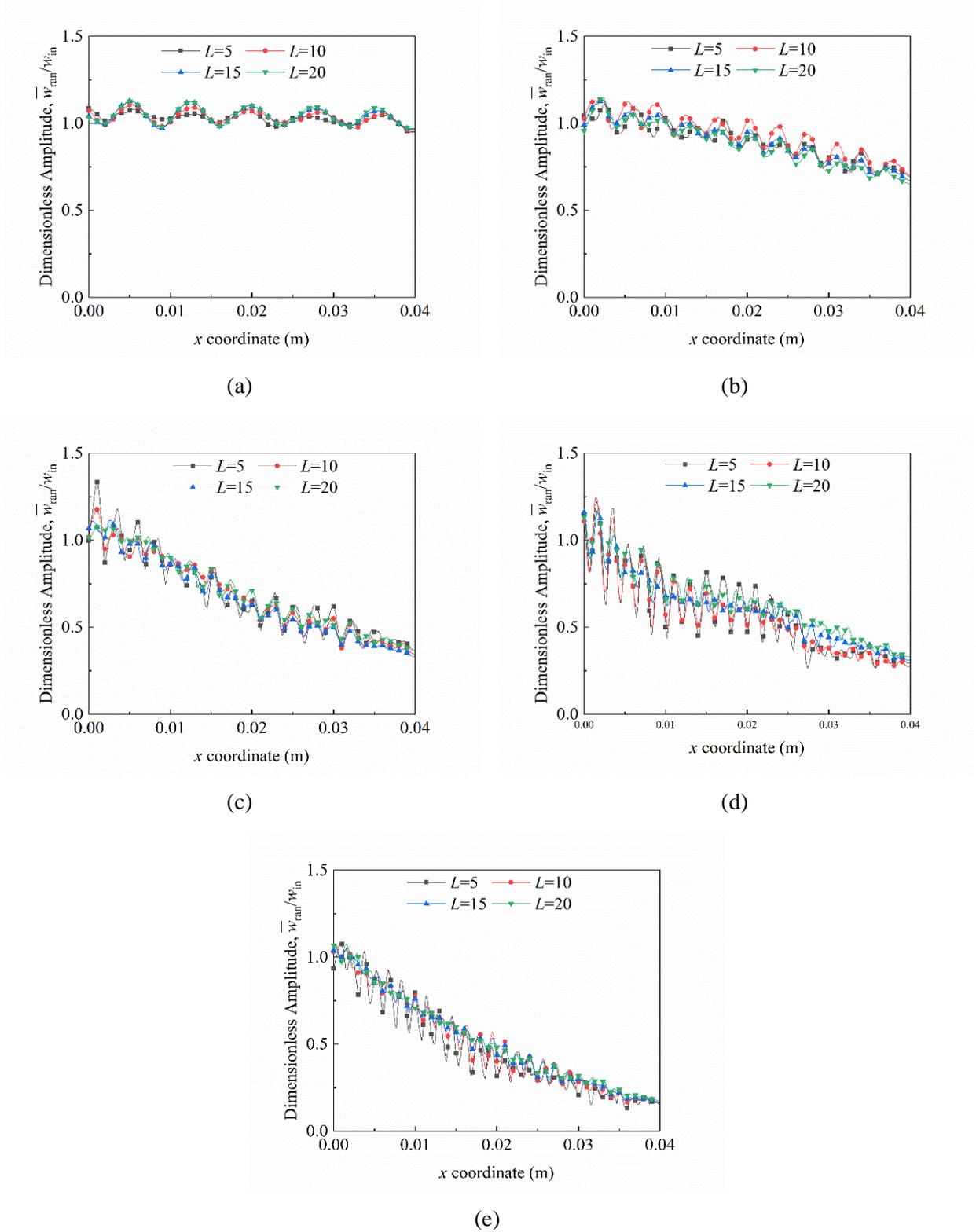

Fig. 8. The dependence of the dimensionless amplitude of the final averaged displacement curves along $x$- coordinate on the number of random cavity distributions ($L$) for different incident wavenumbers, (*a*) $ka$=0.24; (b) $ka$=0.48; (c) $ka$=0.72; (d) $ka$=0.96; (e) $ka$=1.20.

According to Fig. 8, the dimensionless amplitude of the final averaged displacement curve gradually becomes convergent with the increase of $L$. When the dimensionless incident wavenumber $ka$≤1.20, the curves keep nearly the same at *L*=15 and *L*=20. Thus, *L*=20 is selected to calculate the dimensionless amplitude of the final averaged displacement.

The dimensionless amplitude of the final averaged displacement from Eq. (16) with *L*=20 is summarized in Fig. 9.

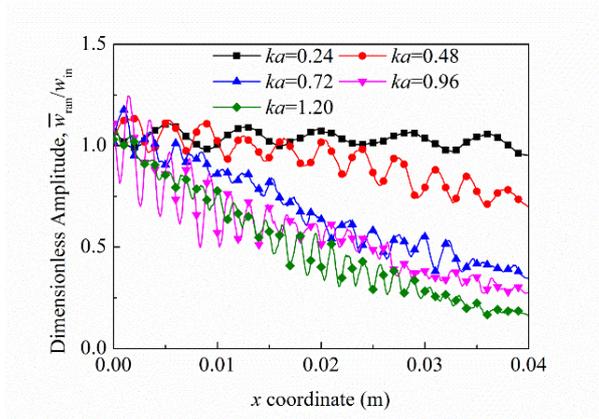

Fig. 9. The dimensionless amplitude of the final averaged displacement along *x*-coordinates for *L*=20 and different incident wavenumbers.

As shown in Fig. 9, all dimensionless amplitude of the final averaged displacement curves along *x*-coordinate for different dimensionless incident wavenumbers start from unity and gradually decrease along *x*-coordinate. For the dimensionless incident wavenumber of 0.24, 0.48, 0.72, 0.96 and 1.20, the dimensionless amplitude of the final averaged displacement at *x*=40 mm are 0.95, 0.70, 0.34, 0.27 and 0.16, respectively, which implies that a larger incident wavenumber or a smaller incident wavelength will lead to a more significant attenuation effect on the amplitude of the final averaged displacement in the direction of wave propagation.

## 4. Finite element verification

The effectiveness of the proposed calculation method in Section 2.1 and the results in Section 3 on quantifying the attenuation effects in the model with cylindrical cavities are verified using FEM with Abaqus.

A 2-dimensional model with cylindrical cavities is built for finite element verification, as shown in Fig. 10. 50 cavities are embedded inside the model with a distribution generated using RSA method. All cavities have the same radius of 0.6 mm as that used in Section 2.1.

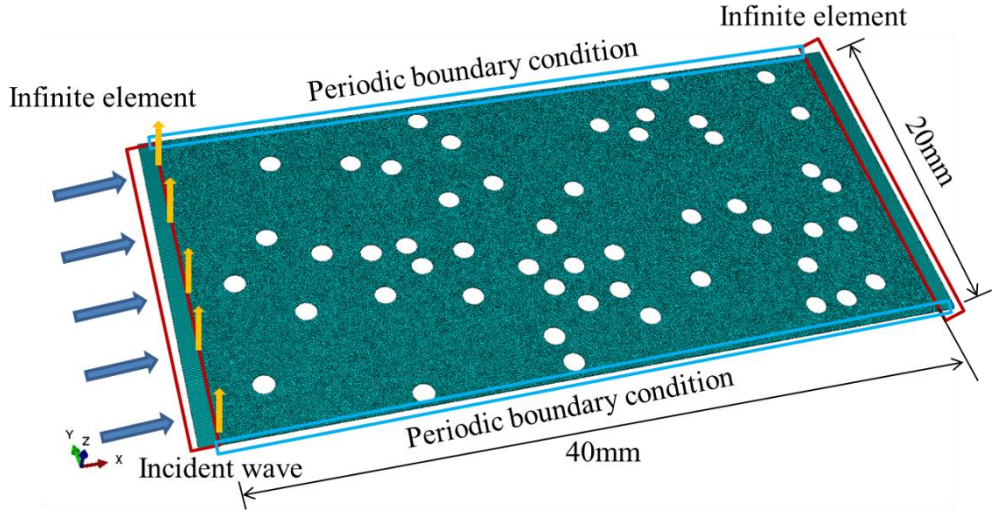

Fig. 10. FEM with random cavities.

The FEM with random cylindrical cavities has length of 40 mm, width of 20 mm and thickness of 0.1 mm. The size of elements is globally selected as 0.1mm × 0.1mm × 0.1mm based on a mesh sensitivity convergence study to balance the computational time and accuracy. There are 88783 elements in total in the FEM. All element motions are allowed only in $z$ direction (i.e. $u=v=0$ and $w \neq 0$ in $x$, $y$ and $z$ directions, respectively) to simulate the SH wave propagation. Periodic boundary condition is added to the upper and lower boundaries to simulate the infinite repetition of the model in $y$ direction.

The FEM of the representative segment is joined by infinite elements (CIN3D8) at its two ends in order to apply non-reflective boundary conditions. The load is applied on the left boundary between the finite elements and infinite elements to generate the incident wave as shown in Figs. 11 (a). The loads applied on the model are concentrated forces on every joining node between finite and infinite elements of each boundary element, as shown in Fig. 11 (b).

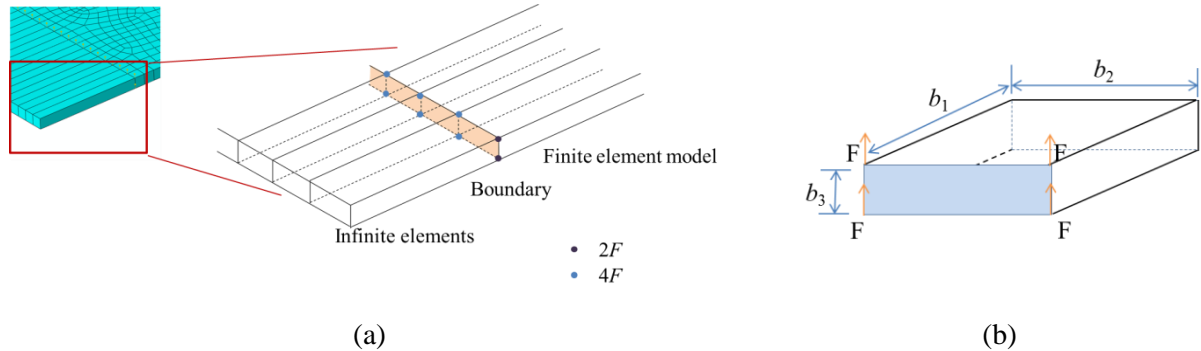

(a)            (b)

Fig. 11. The load applied on the FEM: (a) load position; (b) loads on an element.

The relationship between the load value and the amplitude of the incident wave is calculated as

$$F = \frac{1}{4} \cdot \mu \cdot k \cdot w_{\text{in}} \cdot b_2 \cdot b_3 \tag{19}$$

where $F$ is the value of the load applied on an independent element; $\mu$ is the shear modulus of the element material; $w_{\text{in}}$ is the amplitude of the incident displacement wave; $b_1$, $b_2$, $b_3$ are respectively the length, width and height of the element. The detailed derivation of Eq. (19) is illustrated in Appendix-1.

In this study, the amplitude of the incident wave $w_{in}$ is set as unity. As there is only one element in the thickness direction, the nodes on upper and lower boundaries in *y* direction at two ends (*x* direction) connect two adjacent elements. Thus, the loads applied on the end nodes on upper and lower boundaries (y direction) equal 2*F*. The rest end nodes are all shared by four elements, and therefore, each node is subjected to a load of 4*F*.

The steady state dynamic behaviour of the time-harmonic SH wave propagating in the whole model can be calculated using 'steady-state dynamics, direct' step procedure in Abaqus. The displacement wave field at equally-spaced 17 sections from the left end to the right end are directly extracted by post-processing in Abaqus.

The dimensionless amplitude of the sectional displacements are calculated using the FEM in Abaqus with the same given cavity distribution as those used in MSM in Section 2. In the range of the dimensionless incident wavenumber *ka* from 0.24 to 1.20, the comparison of the dimensionless amplitude of the sectional displacement between the MSM predictions from Eq. (14) and the FEM results is shown in Fig. 12.

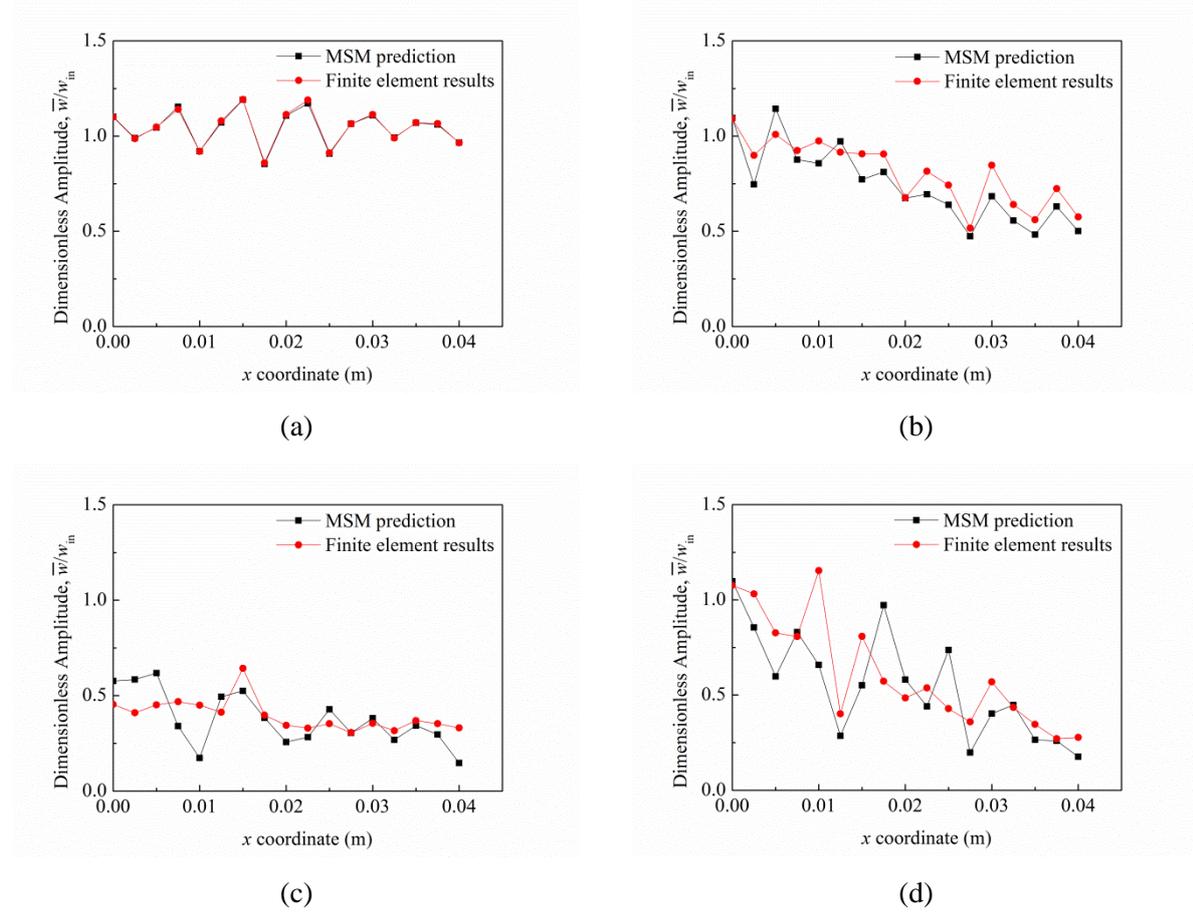

(a)  (b)  (c)  (d)

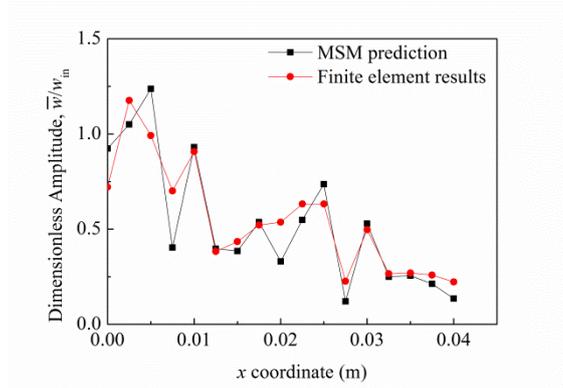

(e)

Fig. 12. The comparison of the dimensionless amplitude of the sectional displacement between the MSM prediction and FEM results and different dimensionless incident wavenumbers of (a) *ka*=0.24; (b) *ka*=0.48; (c) *ka*=0.72; (d) *ka*=0.96; (e) *ka*=1.20.

As shown in Fig. 12, the dimensionless amplitude of the sectional displacement curves predicted by MSM and FEM have very close attenuation along *x*-coordinate, indicating the validity of the collaboration method used in Section 2 when the dimensionless incident wavenumber *ka* ranges from 0.24 to 1.20. All dimensionless amplitudes of the sectional displacement curves show a trend of gradual decrease along the wave propagation in *x*-direction. The initial values of the dimensionless amplitude of the sectional displacement at *x*=0 are all around 1.0 when the dimensionless incident wavenumbers are 0.24, 0.48, 0.96 and 1.20. However, when the dimensionless incident wavenumber *ka*=0.72, the dimensionless amplitude of the sectional displacement at *x*=0 is around 0.5, which is obviously smaller than the corresponding values of the rest curves. This sudden drop of the initial value may be attributed to the use of periodical boundary condition for the whole model. This phenomenon would not happen in a model with cavities distributions introduced in full randomness forthye same model, which is impossible to be implemented in a real simulation, and therefore, an average result from Monte-Carlo study with 20 random distributions is needed.

**5. Homogenisation**

5.1. Effective material properties

The dispersion and attenuation effects are influenced by a variety of independent parameters including the radius and average distance of the embedded cavities, the wavelength, the frequency and amplitude of the incident wave, the Young's modulus, density and Poisson's ratio of the matrix. These parameters have been regrouped into non-dimensional numbers using dimensional analysis (Appendix 2) to reduce the dimension of the parametric space. A large scale parametric analysis is necessary to understand the influence of the meso-scale factors and material properties on the macroscopic wave propagation behaviour, which may benefit from the reduced dimension of the dominant parametric space. However, a general parametric analysis is outside the scope of this study. Instead, *ka* has been selected as a single parameter to show how a homogenized model with damping is obtained to describe the wave propagation behaviour in a model embedded with cavities.

As shown in Section 3, the dispersion effect is negligible when the dimensionless incident wavenumber $ka$ ranges from 0.24 to 1.20 while attenuation effect gradually increases with the increase of the dimensionless incident wavenumber.

To describe the oscillating and gradually decaying wave in macroscopic scale, the homogenized (averaged) displacement defined in Section 2 and used in Sections 3 and 4 is described by

$$\frac{\overline{w}_{ran}}{w_{in}} = \alpha e^{-Ax} \tag{20}$$

where $\alpha$ is a correction coefficient since the dimensionless amplitude of the final displacement $\overline{w}_{ran}/w_{in}$ is not exactly unity at x=0, and

$$A = A_1 i + A_2 \tag{21}$$

is a complex number, in which $A_1 = k_{eff}$ is the imaginary part of $A$ determining the dispersion effect and $A_2$ is the real part of $A$ determining the attenuation effect. Both $A_1$ and $A_2$ can be generally considered as functions of dominant non-dimensional numbers derived in Appendix 2. In this study when other parameters are fixed (see Appendix 2), it has been shown in Section 3.1 that $A_1 = k_{eff}$ is a constant, i.e. the dispersion effect is negligible. $A_2$ as a function of $ka$ could be obtained from fitting the data set of amplitude and propagation distance by $\frac{\overline{w}_{ran}}{w_{in}} = e^{-A_2 x}$ based on least squares method, as shown in Fig. 13.

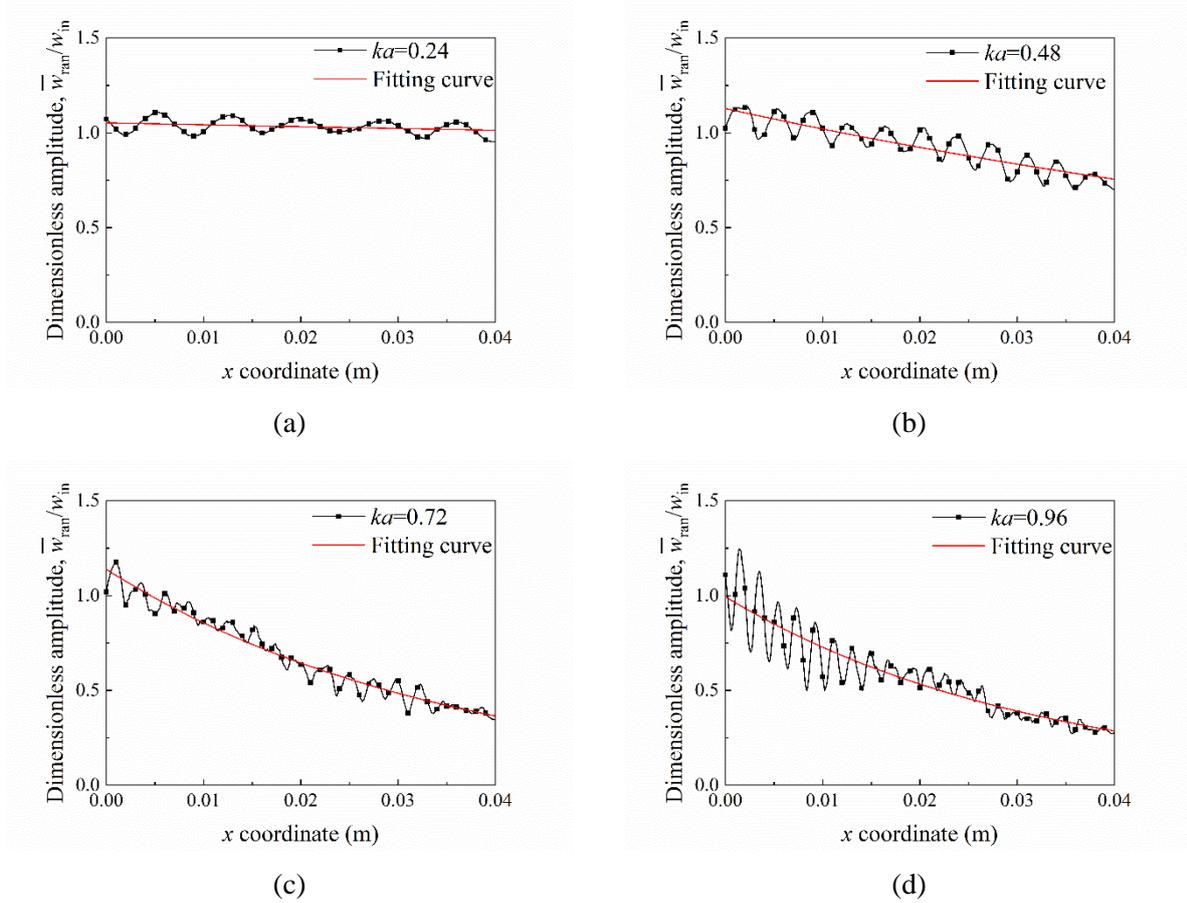

(a)  (b)

(c)  (d)

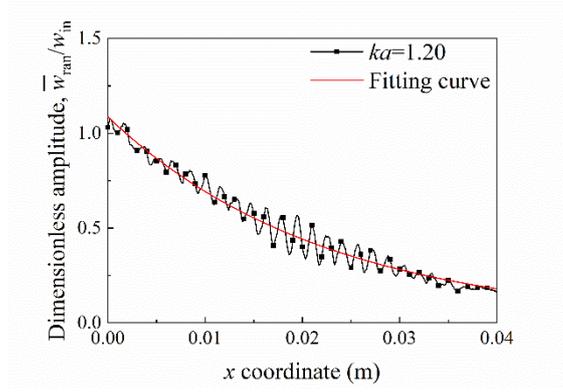

(e)

Fig. 13. The fitting results to calculate the attenuation coefficient $A_2$ for various dimensionless incident wavenumbers, (a) $ka$=0.24; (b) $ka$=0.48; (c) $ka$=0.72; (d) $ka$=0.96; (e) $ka$ =1.20.

Recording the fitting results at different dimensionless incident wavenumbers, the change of attenuation coefficient $A_2$ with the dimensionless incident wavenumber is shown in Fig. 14, where the attenuation coefficient generally shows a linear increase with the incident wavenumber when the dimensionless incident wavenumber $ka$ is less than 1.20. The slope of the fitting line is 34.56.

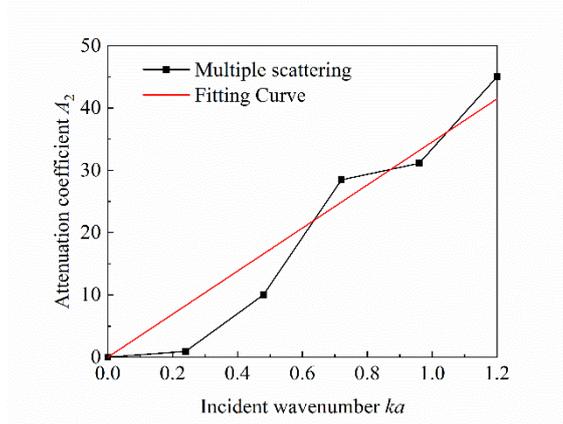

Fig. 14. The change of the attenuation coefficient with dimensionless incident wave wavenumber.

The above-described attenuation effects caused by wave scattering in a porous material can be represented by a structural damping model in a homogenized macroscopic elastic medium.

Based on the assumption that the amplitude of the damping force is proportional to the force caused by elastic deformation and the direction of the damping forces is opposed to the velocity, the dynamic equation is given by (Abaqus, 2014)

$$M\ddot{U} + KU + i \cdot s \cdot KU = 0 \qquad (22a)$$

or

$$M\ddot{U} + (1 + i \cdot s) \cdot KU = 0 \qquad (22b)$$

where $U$ is displacement; $M$ is the mass matrix; $K$ is the stiffness matrix; $s$ is the structural damping coefficient; $i$ is the imaginary unit.

The solution of Eq. (22) has the form of

$$U = e^{i \cdot (\omega t - k_{mac} \cdot x)} = e^{-k_{mac} \cdot x} \cdot e^{i\omega t} \tag{23}$$

where the macroscopic wavenumber is

$$k_{mac} = k_{eff}\sqrt{1 + i \cdot s} = k_{eff}(B_1 i + B_2) \tag{24}$$

for wave propagation in damped elastic material. $k_{mac}$ is a complex number, and $k_{eff}B_1$ and $k_{eff}B_2$ are respectively its imaginary and real parts. $s$ can be expressed explicitly as

$$s = 2B_1 B_2 + (B_1^2 + 1 - B_2^2)i \tag{25}$$

The first term of Eq.(23) has the same form of Eq. (20), and Eq.(24) can be corresponded to Eq.(21). Thus, Eq. (25) can be used to describe the homogenized dynamic behaviour of the waves propagating in the model in Section 3, and the fitting results for Eq. (21) based on Fig. 14 can be used to determine $B_1$ and $B_2$, i.e.

$$k_{eff}B_1 = A_1 \tag{26a}$$

$$k_{eff}B_2 = A_2 \tag{26b}$$

Then, $B_1$ and $B_2$ could be calculated from Eq. (26a, b) as $B_1$=0.0207 and $B_2$=1.000. $s$ is determined by Eq.(25) to be 0.041+0.0002i. Then, the imaginary part of $s$ is neglected and $s$ is selected as 0.041. It can be concluded that the meso-scale wave scattering phenomenon can be represented by a macroscopic structural damping in an equivalent homogenized model.

The volume porosity of the homogenized model of a porous medium with cylindrical cavity of radius a is

$$\eta = N \cdot \frac{\pi \cdot a^2}{H \cdot T} \tag{27}$$

where $N$ is the total number of cavities; $H$ and $T$ are the height and length of the porous medium, respectively.

The effective density of the homogenized model is

$$\rho_{eff} = \rho(1 - \eta) \tag{28}$$

where $\rho$=2.7×10³ kg/m³ is the material density of the solid phase of the porous medium. For the presently studied problem, $a$=0.6mm, $H$=20mm, $T$=40mm and $N$=50, thus, the volume porosity $\eta$=0.071 and effective density $\rho_{eff}$=2.5×10³ kg/m³.

The dynamic effective shear modulus of the porous media can be determined by

$$\mu_{eff} = c_{eff}^2 \cdot \rho_{eff} \tag{29}$$

where the effective wave speed $c_{eff}$ can be calculated from Eq.(18) for a given incident wave frequency, in which $k_{eff}$ is calculated as the average of wavenumbers measured between neighboring peaks on $\bar{w}^r/w_{in}$ –$x$ graph. Since the wave dispersion is negligible in the present problem, $k_{eff}/k$ has little change with the incident wave frequency (wavenumber). Therefore, $c_{eff}$ is almost independent of the incident wave frequency. Herein, an average $c_{eff}$ is calculated using the average of $k_{eff}$ for all dimensionless wavenumbers considered in this study. The dynamic effective shear modulus is determined as 24.05 GPa. A method based on effective cross-sectional area is also used to determine the effective shear modulus, as shown in Appendix 3, and a value of $\mu_{eff}$ = 24.5 GPa is obtained, which is very close to the value of the dynamic effective shear modulus. This implies that the effective Young's modulus and density measured for a porous material can be used to determine wave speed when porosity-induced wave dispersion is negligible.

The effective Poisson's ratio $v$ of the porous medium is assumed to be the same as that of the matrix, namely 0.3. Then the dynamic effective Young's modulus can be calculated as

$$E_{eff} = 2(1+v)\mu_{eff} \tag{30}$$

Therefore, the homogenized (SH wave) model of the porous medium is an elastic model with effective density 2.5×10³ kg/m³, the effective Young's modulus 62.52 GPa, Poisson's ratio 0.3 and with structural damping 0.041 based on the calculations in Eq. (26-30).

5.2. Homogenized FE model

Structural damping modelling is also available in Abaqus for structural dynamics analysis. To check the effectiveness and applicability of the homogenized model based on effective material properties, a homogeneous (SH wave) FEM corresponding to the meso-scale FEM presented in Section 4 is established to simulate the wave propagation in an unbounded material (in $x$-direction) with the consideration of structural damping, as shown in Fig. 15. The model has the same geometrical dimensions as those in the model used in Sections 2.2 and 4. The FEM of the representative segment is sandwiched between infinite elements in $x$ direction and limited to move in $z$ direction. Then, the density, Young's modulus and Poisson's ratio of the homogeneous FEM are respectively set as 2.5×10³ kg/m³, 62.52 GPa and 0.3. The structural damping of the homogeneous FEM is set as 0.041.

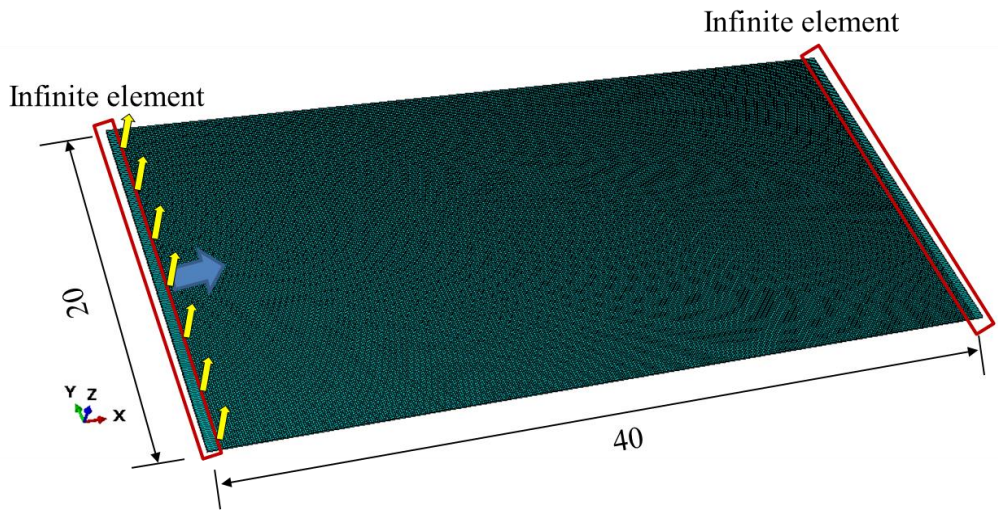

Fig. 15. The FEM to simulate the macroscopic wave propagation in a continuum with damping

The size of elements in the model is set as 0.1mm×0.1mm×0.1mm. There are 80000 elements in the middle segments, all of which are of element type C3D8R. The analyses are conducted in steady-state dynamics, direct step. A time-harmonic SH wave with unity amplitude is generated from the boundary between the left end of the model and the outside infinite elements, as shown in Fig.15. The origin of the coordinate system is set at the centre of the FEM..

5.3. Comparison between MSM and FEM

The dimensionless final averaged displacements at different sections in the porous media calculated from Eq. (15a) of MSM with L=20 are compared with the real component of displacements at 400 equally-spaced sections that extracted from the homogeneous elastic FEM with structural damping, as shown in Fig. 16.

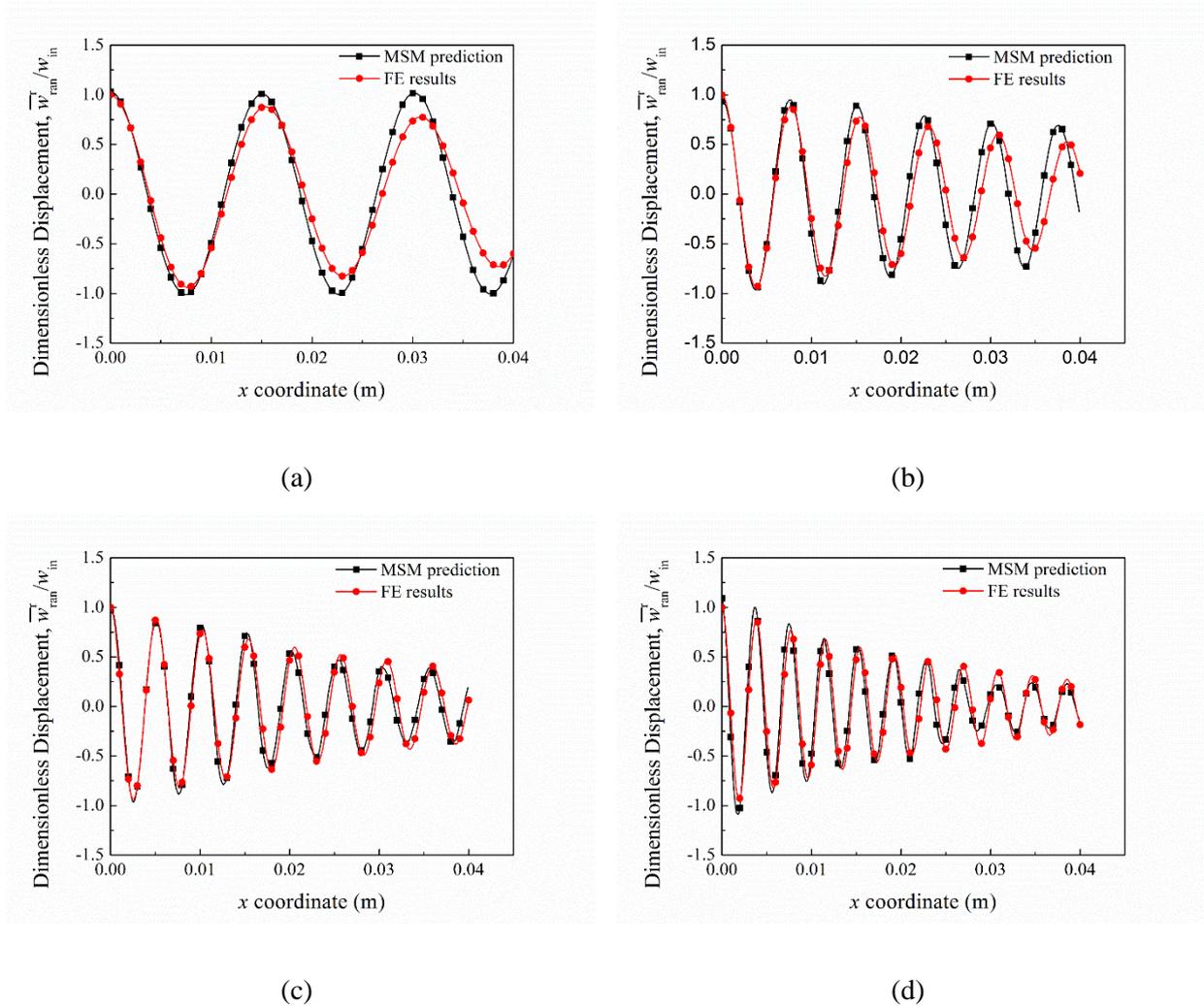

(a)  (b)

(c)  (d)

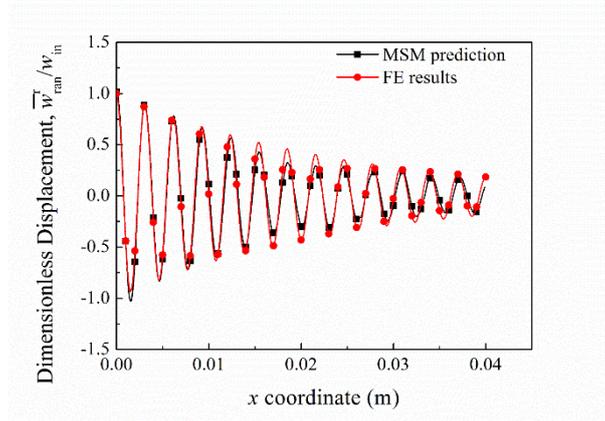

(e)

Fig. 16. Comparison of the dimensionless final averaged displacements from MSM and FEM. (a) *ka*=0.24; (b) *ka*=0.48; (c) *ka*=0.72; (d) *ka*=0.96; (e) *ka*=1.20.

As shown in Fig. 16, all dimensionless final averaged displacement curves along *x*-coordinate from FEM shows a similar trend with those from MSM. The position of peak points in Fig. 16(a) shows obvious difference between the two curves than those in Fig. (16b, c, d and e), which can be attributed to a relatively larger error of the dimensionless wavenumber $k_{\text{eff}}/k$ to the average value adopted in the homogenized model at *ka*=0.24 in Table 1. Similarly, the larger difference of the attenuation coefficient between the calculated result and the fitting result in Fig. 14 makes a larger difference of the amplitudes of peak points between MSM prediction curves and FEM curves in Fig. 16(a).

The fitting quality of the FEM curves in Fig. 16 can be quantified by coefficients of determination ($R^2$), as shown in Table 2.

Table 2. The coefficients of determination ($R^2$) for the FEM results and MSM predictions.

| Dimensionless wavenumber *ka* | 0.24 | 0.48 | 0.72 | 0.96 | 1.20 |
|---|---|---|---|---|---|
| Coefficients of determination $R^2$ | 0.95 | 0.85 | 0.94 | 0.89 | 0.96 |

As shown in Table 2, when the dimensionless incident wavenumber *ka* ranges from 0.24 to 1.20, the coefficient of determination for the dimensionless final averaged displacements predicted by FEM and MSM are respectively 0.95, 0.85, 0.94, 0.89 and 0.96, which shows high effectiveness of the FEM structural damping model to describe the wave propagation behaviour in porous media.

## 6. Conclusions

This paper combines eigenfunction method and collaboration method to give an analytical prediction of the steady-state wave field in a porous medium excited by a time-harmonic elastic incident wave, which is verified using a meso-scale finite element analysis. A Monte Carlo study is then conducted to get an average wave propagation results in a model with randomly-distributed cavities, relating to the homogenised macroscopic wave propagation behaviours. Finally, homogenized wave propagation behaviour is demonstrated by an equivalent macroscopic structural damping model. The main findings are

(i) The sectional averaged SH wave propagation behaviour in porous media with sufficiently large number of random cavities can be analytically predicted by MSM to represent the macroscopically homogenised SH wave behaviour.

(ii) The effectiveness of MSM for the study of multiple wave reflections by meso-scale inclusions can be verified by meso-scale finite element modelling.

(iii) The homogenized steady-state SH wave behaviours in porous media can be effectively described by a homogenous elastic model with effective material properties and a pseudo structural damping.

The presented methodology may be further developed to understand the dominant meso-scale influential factors and their effects on the macroscopic material behaviours and wave propagations for more complex porous and heterogeneous media.

**Acknowledgemen**t: The first author acknowledges the support from the President's Doctoral Scholarship (PDS) awarded by The University of Manchester.

**Appendix 1 The calculation of applied force to replace the incident wave**

The shear stress is proportional to the shear strain,

$$\tau = \mu\gamma \qquad (A.1)$$

where $\tau$ is the shear stress, $\gamma$ is the shear strain, $\mu$ is the shear modulus of the matrix. An incremental expression of Eq. (A.1) is

$$d\tau = \mu d\gamma \qquad (A.2)$$

Consider the case that an incident wave propagates from infinite distance in a 1D model. A small segment BC in the wave propagation path is selected, as shown in Fig. A.1. Herein $\lambda$ is the incident wavelength; $F_B$ and $F_C$ are respectively the shear forces on B surface and C surface; $b_1$ is the distance between B and C along the wave propagation direction; $h$ is the displacement difference in the out-off-plane direction; $\Delta\theta$ is the angle of inclination between BC and the wave propagation direction; $\phi$ is the phase difference between B and C; $\varphi$ is the phase difference between B and the next point of zero displacement.

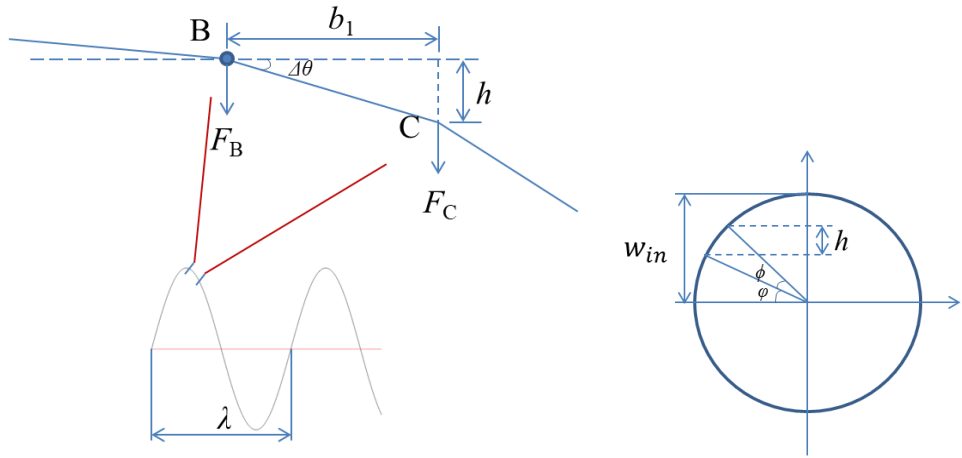

Fig. A.1. The calculation of the equivalent load boundary condition to displace boundary condition.

The 1D simulation of shear wave propagation is performed using Abaqus by constraining the motion of elements in z direction, as shown in Fig. 11(b). Same load F is applied to all nodes on surface B.

If the element is small enough,

$$d\tau = \frac{\Delta F}{b_2 \cdot b_3} \qquad (A.3a)$$

$$d\gamma = \Delta\theta = \frac{h}{b_1} \qquad (A.3b)$$

Then, Eq. (A.2) can be represented by

$$\frac{\Delta F}{b_2 \cdot b_3} = \mu \cdot \frac{h}{b_1} \qquad (A.4)$$

where $b_2$ and $b_3$ are defined in Fig.11(b), and

$$\Delta F = F_\text{B} - F_\text{C} \tag{A.5}$$

The phase difference between B and C is

$$\phi = \frac{b_1}{\lambda} \cdot 2\pi = k b_1 \tag{A.6}$$

where $k$ is the incident wavenumber. Then,

$$h = w_\text{in} \cdot (\sin(\varphi + \phi) - \sin(\varphi)) \tag{A.7}$$

where $w_\text{in}$ is the amplitude of the incident wave.

When the load at B reaches the maximum, the phase difference between B and its nearest zero load surface is π/2. Then, the maximum force $F_\text{p}$ at B could be obtained through integration, i.e.

$$\begin{aligned}\frac{F_P}{b_2 \cdot b_3} &= \mu \cdot \frac{w_\text{in}}{b_1} \cdot \int_0^{\frac{\pi}{2}} [\sin(\varphi + \phi) - \sin(\varphi)]\, d\varphi \\ &= \mu \cdot \frac{w_\text{in}}{b_1} \cdot (\cos(kb_1) - 1 + \sin(kb_1))\end{aligned} \tag{A.8}$$

If $kb_1 \to 0$,

$$\cos(kb_1) - 1 + \sin(kb_1) \approx kb_1 \tag{A.9}$$

thus,

$$\frac{F_P}{b_2 \cdot b_3} \approx \mu \cdot \frac{w_\text{in}}{b_1} \cdot kb_1 = \mu \cdot k \cdot w_\text{in} \tag{A.10}$$

A concentrated force of $F = F_P/4$ is applied on every single node of the boundary elements to generate an incident wave with amplitude of $w_\text{in}$.

**Appendix 2 Dimensional analysis and a harmonic SH wave in a porous media**

For a porous medium in the *x-y* plane with randomly-distributed circular voids characterised by its radius *a* and averaged centre distance *l*, consider the propagation of the *y*-direction averaged SH displacement wave, *w=w(x, t)* in *x*-direction in the porous medium, as shown in Fig. 1. The solid phase of the porous medium is assumed to be isotropic and linear elastic, which can be represented by Young's modulus *E* or the shear modulus $\mu$, density $\rho$ and Poisson's ratio *v*. The porosity of the porous medium is given by

$$\eta = \pi(\frac{a}{l})^2 \tag{A.11}$$

If a harmonic displacement excitation in z-direction ($w_0$) is applied at *x*=0, i.e.

$$w_0 = w_{in}e^{-i\omega t} \tag{A.12}$$

where $w_{in}$ and $\omega$ are the amplitude and angular frequency of the incident wave,

The propagation of the SH wave (*w*) can be completely determined by $w_{in}$, $\omega$, $\mu$, $\rho$, *v*, *a* and *l* in addition to *x* and *t*, i.e.

$$w = w(a, l, w_{in}, v, \omega, \mu, \rho, x, t) \tag{A.13a}$$

or

$$w = w(a, l, w_{in}, v, k_0, \mu, \rho, x, t) \tag{A.13b}$$

when $\omega$ is replaced by *k* based on following relationships

$$k = \frac{2\pi}{\lambda} \tag{A.14}$$

and

$$c = \frac{\omega}{k} = \sqrt{\frac{\mu}{\rho}} \tag{A.15}$$

where $\lambda$, $k$ and $c$ are the wavelength, wavenumber and phase velocity of the harmonic SH wave propagating in the continuum of matrix material. Eqs. (A.14) and (A.15) are applicable for the propagation of harmonic wave in the equivalent continuum medium homogenised from the corresponding porous medium where a subscript 'eff' will be added for the equivalent continuum.

Using dimensional analysis, Eq. (A. 13b) with ten physical quantities can be reduced to following non-dimensional form with seven non-dimensional physical quantities, i.e.

$$\frac{w}{w_{in}} = F(ka, kl, kw_{in}, v, kx, kt\sqrt{\mu/\rho}) \tag{A.16a}$$

or

$$\frac{w}{w_{in}} = F(ka, \eta, \frac{w_{in}}{a}, v, \xi, \tau) \tag{A.16b}$$

when $kl$ is replaced by $\eta$, $kw_{in}$ is replace by $\frac{w_{in}}{a}$, and $\xi = k_0 x$, $\tau = k_0 t\sqrt{\mu/\rho}$.

It is assumed that the propagation of the averaged wave $w=w(x, t)$ in the original porous medium can be described by its corresponding effective wave propagation in the equivalent continuum. Therefore, the effective harmonic SH wave, corresponding to $w=w(x, t)$, in the equivalent medium can be described by

$$w_{\text{eff}} = w_{in} e^{i(k_{\text{eff}} x - \omega t)} \tag{A.17a}$$

or

$$\frac{w_{\text{eff}}}{w_{in}} = e^{i(k_{\text{eff}} x - k_{\text{eff}} c_{\text{eff}} t)} \tag{A.17b}$$

where $k_{\text{eff}}$ is the effective wavenumber and $c_{\text{eff}}$ is the effective wave phase velocity in the equivalent continuum.

In this research, the amplitude of incident wave number is set as unity; The radius of cavities (a) is given; the porosity of the model $\eta$ is fixed as the number of cavities are fixed as 50; the shear modulus $\mu$ and the density $\rho$ are also fixed. Then this study only uses $ka$ as an independent non-dimensional number for the analysis of the effective wave propagation.

**Appendix 3 The calculation of effective shear modulus**

An equivalent self-consistent model (SCM) of the studied porous medium in Section 3 is represented by a cube with a central cylindrical void, shown in Fig. A.2.

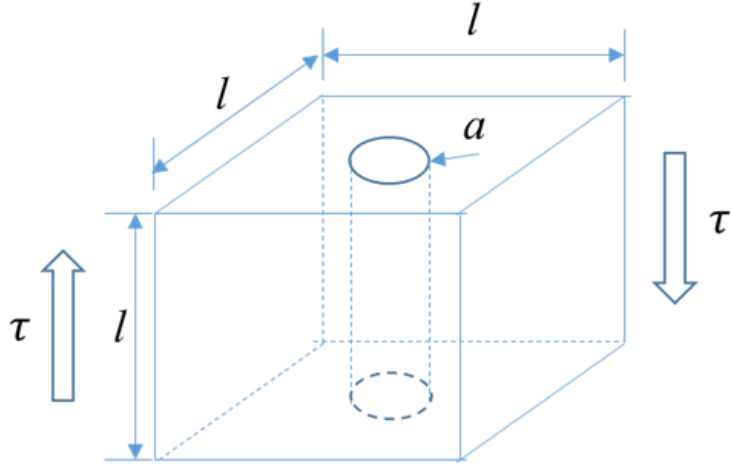

Fig. A.2. The 3D self-consistent model of the studied porous medium.

The length, height and width of the SCM are $l$ and the radius of the embedded cavity is $a$. The sectional area perpendicular to the SH wave propagation direction is

$$A = l^2 \tag{A.18}$$

The volume porosity of the SCM is the same as the porosity of the porous medium in Section 3 calculated by Eq. (27) as 0.071, i.e.

$$\eta = \frac{\pi a^2}{l^2} = 0.071 \tag{A.19}$$

Then, the ratio of the radius of the embedded cavity $a$ to the length of the SCM $l$ is obtained as

$$\frac{a}{l} = 0.15 \tag{A.20}$$

The effective shear modulus of the SCM, $\mu_{\text{eff}}$, can be calculated through the total shear displacement $\delta$ in the wave propagation direction as follows.

$$\mu_{\text{eff}} = \frac{\tau}{\gamma} = \frac{\tau l}{\delta} \tag{A.21}$$

A Cartesian coordinate system is established at the centre of the cross section in wave propagation direction, as shown in Fig. A3.

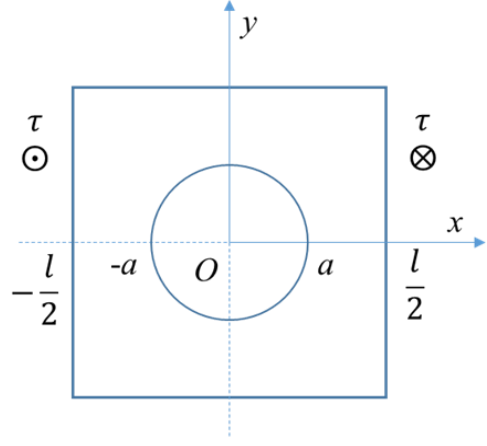

Fig. A.3. Top view of the cross-sectional shape of the SCM and the coordination system

The total shear displacement in the SCM is calculated as

$$\delta = \frac{\tau}{\mu}\left[-r - \left(-\frac{l}{2}\right)\right] + \int_{-r}^{r}\frac{\tau(x)}{\mu}dx + \frac{\tau}{\mu}\left(\frac{l}{2} - r\right)$$
$$= \frac{2\tau}{\mu}\left(\frac{l}{2} - r\right) + \int_{-r}^{r}\frac{\tau(x)}{\mu}dx \quad (A.22)$$

where $\tau$ is the external shear stress applied on the SCM, $\mu$ is the shear modulus of matrix in the SCM. $\tau(x)$ is the shear stress distribution at $-a \leq x \leq a$. The equilibrium of shear across different $x$ coordinate makes

$$\tau \cdot \frac{l}{2} = \tau(x)\left(\frac{l}{2} - y\right), \text{at} -a \leq x \leq a \quad (A.23)$$

The coordinates of the circular void between $-a \leq x \leq a$ are

$$x = r \cdot \cos(\theta), \quad (A.24a)$$

$$y = r \cdot \sin(\theta) \quad (A.24b)$$

Then, Eq. (A. 22) can be rewritten as

$$\delta = \frac{2\tau l}{\mu}\left(\frac{1}{2} - \frac{a}{l}\right) + \frac{\tau l}{\mu}\int_{\pi}^{0}\frac{-(\frac{a}{l})\sin(\theta)}{1 - \frac{2a}{l}\cdot\sin(\theta)}d\theta \quad (A.25)$$

Substitute Eq. (A.20) into Eq. (A.25),

$$\delta = 1.10 \cdot \frac{\tau l}{\mu} \quad (A.26)$$

Substitute Eq. (A.26) into Eq. (A.21),

$$\mu_{\text{eff}} = \frac{\mu}{1.10} = 24.5 \text{ GPa} \quad (A.26)$$